\documentclass[12pt]{amsart}

\usepackage{geometry}                
\usepackage{graphicx}
\usepackage{amssymb}

\newtheorem{theorem}{Theorem}[section]

\newtheorem{lemma}[theorem]{Lemma}
\newtheorem{claim}[theorem]{Claim}
\newtheorem{fact}[theorem]{Fact}
\newtheorem{corollary}[theorem]{Corollary}
\newtheorem{proposition}[theorem]{Proposition}

\theoremstyle{definition}

\newtheorem{problem}[theorem]{Problem}
\newtheorem{example}[theorem]{Example}

\theoremstyle{remark}
\newtheorem{remark}[theorem]{Remark}

\numberwithin{equation}{section}

\newcommand{\cH}{\mathcal{H}}

\newcommand{\cA}{\mathcal{A}}

\newcommand{\cB}{\mathcal{B}}

\newcommand{\D}{\mathbb{D}}

\newcommand{\C}{\mathbb{C}}
\newcommand{\T}{\mathbb{T}}

\newcommand{\N}{\mathbb{N}}

\newcommand{\R}{\mathbb{R}}

\newcommand{\Om}{\Omega }

\newcommand{\diam}{\operatorname{Diam}}
\newcommand{\ndiam}{\operatorname{n-Diam}}

\newcommand{\rea}{\operatorname{Re}}
\newcommand{\dist}{\operatorname{dist}}
\newcommand{\length}{\operatorname{Length}}
\newcommand{\ima}{\operatorname{Im}}

\newcommand{\area}{\operatorname{Area}}

\newcommand{\bd}[1]{\partial #1}

\newcommand{\capa}{\operatorname{Cap}}
\newcommand{\rad}{\operatorname{Rad}}

\newcommand{\defeq}{\mathrel{\mathop:}=}

\title{Area, capacity and diameter versions of Schwarz's Lemma}
\author[Burckel]{Robert B. Burckel}
\address{Department of Mathematics, Cardwell Hall, Kansas State University,
Manhattan, KS 66506, USA}
\email{burckel@math.ksu.edu}

\author[Marshall]{Donald E. Marshall}
\address{Department of Mathematics, Box 354350
University of Washington 
Seattle, WA 98195-4350 USA}
\email{marshall@math.washington.edu}

\author[Minda]{David Minda}
\address{Department of Mathematical Sciences, 
 University of Cincinnati, 
 PO Box 210025,
 Cincinnati, OH 45221-0025 }
\email{david.minda@uc.edu}

\author[Poggi-Corradini]{Pietro Poggi-Corradini}
\address{Department of Mathematics, Cardwell Hall, Kansas State University,
Manhattan, KS 66506, USA}
\email{pietro@math.ksu.edu}

\author[Ransford]{Thomas J.~Ransford}
\address{D\'epartement de math\'ematiques et de statistique,
  Universit\'e Laval, Qu\'ebec (QC), G1K 7P4, Canada}  
\email{thomas.ransford@mat.ulaval.ca}

\thanks{The authors thank
    A.~F.~Beardon for his comments and suggestions on an earlier version.}
\subjclass[2000]{30C80}
\date{July 27, 2007}     
\begin{document}
\begin{abstract}
  The now canonical proof of Schwarz's Lemma appeared in a 1907
  paper of Carath\'eodory, who attributed it to Erhard
  Schmidt. Since then, Schwarz's Lemma has acquired considerable
  fame, with multiple extensions and generalizations. Much less known
  is that, in the same year 1907, Landau and Toeplitz obtained a
  similar result where the diameter of the image set takes over the
  role of the maximum modulus of the function. We give a new
  proof of this result and extend it to include bounds on the
  growth of the maximum modulus. We also 
develop a more general 
approach in which the size of the image is estimated in several
geometric ways via notions of radius, diameter, perimeter, area,
  capacity, etc...
\end{abstract}

\maketitle
\baselineskip=18pt

\section{Introduction}
\subsection{Schwarz's Lemma}
First, let us set the following standard notations: $\C$ denotes the
complex numbers, $\D\defeq\{z\in \C: |z|<1\}$ is the open unit disk, and
$\T\defeq\{z\in \C: |z|=1\}$ is the unit circle. Moreover, for $r>0$, we let
$r\D\defeq\{z\in\C: |z|<r\}$ and $r\T\defeq\{z\in \C: |z|=r\}$. Also,
we will say that
a function is {\em linear} if it is of the form $f (z)= az+b$ with $a,b\in
\C$ (in particular, it may be constant).

Given a function $f$ that is analytic in $\D$, and given the exhaustion
$\{r\D \}_{0\leq r \leq 1}$, consider the corresponding image-sets
\[
f (r\D)=\{w\in \C: \mbox{ there is at least one }z\in r\D\mbox{ such
that }f (z)=w\}.
\]
Let us emphasize that we will not consider ``multiplicity'' in this
paper. So $f (r\D)$ denotes a family of open connected sets in $\C$
that are increasing with $r$. The goal is to fix a geometric quantity
so as to measure the size of $f(r\D)$ and study how it varies with
$r$. In particular, it turns out that linear functions seem always to exhibit
a uniquely exceptional behavior.

To illustrate this point of view, we first consider the famous Schwarz's
Lemma. Introduce the following notion of ``radius'':
\begin{equation}\label{eq:max}
\rad f (r\D)\defeq \sup_{|z|<r}|f(z)-f(0)|.
\end{equation}
Geometrically, $\rad f (r\D)$ is the radius of the smallest disk
centered at $f (0)$ which contains  $f (r\D)$.
As a point of warning, the way Schwarz's Lemma will be presented below
might look unusual, but the proof is exactly the same.
\begin{theorem}[Schwarz's Lemma]\label{thm:sl}
Suppose $f$ is analytic on the unit disk $\D$. 
Then the function $\phi_{\rad} (r)\defeq  r^{-1}\rad f (r\D)$ is 
strictly increasing for $0<r<1$, except when $f$ is linear, in which
case $\phi_{\rad}$ is constant. Moreover,
$\lim_{r\downarrow 0}\phi_{\rad} (r) =|f^\prime(0)|.$
\end{theorem}
\begin{corollary}\label{cor:sl}
Suppose $f$ is analytic on the unit disk $\D$ with $\rad f (\D)=1$. Then
\begin{align}
&\label{eq:growth}
\rad f (r\D)\leq r \qquad  \mbox{for every $0<r<1$, and}\\
& \label{eq:der} |f^\prime(0)|\leq 1. 
\end{align}
Moreover, equality holds in 
{\rm (\ref{eq:growth})} for some $0<r<1$, or in {\rm (\ref{eq:der})},
if and only if $f (z)$ is an Euclidean isometry $a+cz$
for some constants  $a\in \C$, $c\in \T$.
\end{corollary}

The standard way to prove Schwarz's Lemma is to factor 
$f(z)-f(0)=zg(z)$, for some analytic function $g$ and then
apply the maximum modulus theorem to $g$ to deduce that
\begin{equation}\label{eq:g}
r^{-1}\rad f (r\D)=\sup_{|z|<r}|g(z)|.
\end{equation}
This argument first appeared in a paper of
Carath\'eodory \cite{caratheodory1907} where the idea is attributed to
E.~Schmidt.  See Remmert, \cite{remmert1991} p.~272-273, and
Lichtenstein, \cite{lichtenstein1919} footnote 427, for historical
accounts. 

\subsection{The theorem of Landau and Toeplitz}

In a 1907 paper, Landau and Toeplitz replaced the radius
(\ref{eq:max}) by the diameter of the image set
\begin{equation}\label{eq:diam2}
\diam f (r\D)\defeq \sup_{z,w\in r\D}|f(z)-f (w)|.
\end{equation}
\begin{theorem}[Landau-Toeplitz \cite{landau-toeplitz1907}]\label{thm:lt}
Suppose $f$ is analytic on the unit disk $\D$ and $\diam f (\D)=2$.
Then  
\begin{align}
\hbox{\,}&\label{eq:diamineq} \diam f (r\D)
\leq 2r\qquad \mbox{for every $0<r<1$, and}\\
&\label{eq:lt} |f^\prime(0)|\leq 1.
\end{align}
Moreover, equality holds in {\rm
(\ref{eq:diamineq})} for some $0<r<1$, or in {\rm (\ref{eq:lt})}, if
and only if 
$f (z)$ is an
Euclidean isometry $a+cz$
for some constants  $a\in \C$, $c\in \T$.
\end{theorem}

\begin{remark}
The main contribution of the Landau-Toeplitz paper is
perhaps its elucidation of the extremal case.
P\'olya and Szeg\H o mention the inequality (\ref{eq:lt}) on p.~151 and
p.~356 of the classic book \cite{polya-szego1972}, and they cite
the paper of Landau and Toeplitz. However, they say nothing about when
equality holds. It is also worth mentioning the proof of Lemma
2.9 in \cite{gehring-hag:1999jcam}, where F.W.~Gehring and K.~Hag
essentially treat the case of equality in Theorem \ref{thm:lt} in the
special case of one-to-one maps,  
using
the Cauchy-Schwarz inequality and the so-called ``isodiametric'' inequality.
\end{remark}
\begin{remark}\label{rem:triangle}
The growth estimate on the diameter
(\ref{eq:diamineq}) should be viewed in analogy with the classical
growth bound (\ref{eq:growth}). 
Notice, however, that Theorem \ref{thm:lt} covers the case when
$f(\D)$ is an equilateral triangle of side-length $2$, which is of
course not contained in a disk of radius $1$; likewise, when $f (\D)$
is contained in the so-called
Reuleaux triangle that is obtained from the equilateral triangle by
joining adjacent vertices by a circular arc having center at the third vertex. 
\end{remark}

We start by giving a new proof of the
Landau-Toeplitz Theorem that can be used to prove more general cases
as well. Later we will show how the original proof of
Landau and Toeplitz can be adapted to some of these more general
cases. The Landau-Toeplitz approach is more direct but seems to
accomplish less.

\subsection{Higher-diameters and log-convexity}
As we already mentioned in Remark \ref{rem:triangle}, the
Landau-Toeplitz result 
generalizes the bounds on $|f^{\prime} (0)|$ that can be deduced from
Schwarz's Lemma. It is therefore natural to ask if there 
are other conditions on 
analytic functions $f:\D\rightarrow \C$, weaker than $\diam f
(\D)\leq 2$, which imply $|f^{\prime} (0)|\leq 1$ with equality if and
only if $f$ is an Euclidean isometry.
Also, it follows from (\ref{eq:g}) and Hadamard's three-circles theorem that
$\phi_{\rad} 
(r)$ is not only strictly increasing (except when $f$ is linear), but
it is also {\sf log-convex}, i.e., it is a 
convex function of $\log r$. In fact, even more is true: its
logarithm is log-convex. 
Thus, a more general question arises: assuming that $f$ is not linear, is 
the function $\phi_{\diam} (r)\defeq (2r)^{-1}\diam f (r\D)$ strictly
increasing and log-convex?

Other geometric quantities may be used to
measure the size of the image of an analytic function. In this paper
we will focus on $n$-diameter, capacity, area and perimeter. In
\cite{polya-szego1951}  P\'olya and Szeg\H o consider also other
quantities such as 
the moment of inertia, the torsional rigidity, and the
principal frequency. Such topics deserve to be explored but we reserve
to do this in another paper.

We focus at first on the
so-called higher-order diameters, which are defined for sets
$E\subset \C$ as 
follows: fix $n=2,3,4,\dots$, then
\[
d_{n} (E)\defeq \sup\left(\prod_{j<k}|\zeta_{j}-\zeta_{k}|
\right)^{\frac{2}{n (n-1)}}, 
\]
where the supremum is taken over all $n$-tuples of points from
$E$. We say $d_n(E)$ is the {\sf $n$-diameter} of $E$. 
Note that $d_{2} (E)=\diam E$, and that $d_{n} (E)$ is weakly
decreasing in $n$. 
Hence $d_{\infty} (E)\defeq \lim_{n\rightarrow \infty}d_{n} (E)$ is
well-defined and is
called the {\sf transfinite diameter} of $E$. 
It turns out that the
transfinite diameter $d_\infty(E)$ coincides with the logarithmic
capacity $\capa(E)$; see the Fekete-Szeg\H o Theorem p.~153 of \cite{R}.
It can also be shown, see
Fact \ref{fact} below, that $d_n(\D)=n^{1/(n-1)}$.

The following inequality is due to P\'olya; see \cite{POL}
or \cite[p.145]{R}.
For a compact set $E$ in $\mathbb{C}$, 
\begin{equation}\label{eq:polya}
\area E \leq \pi\, \capa^2 E.
\end{equation}
Equality holds for a closed disk. 
Moreover, Corollary~6.2.4 of \cite{R}
asserts that
\[
\capa(E)\leq\frac{d_n(E)}{n^{1/(n-1)}}.
\]
Hence, combining with (\ref{eq:polya}), we obtain,
for $n=2,3,...$,
\begin{equation}\label{eq:areadn}
\area E\leq \pi\frac{d_n(E)^2}{n^{2/(n-1)}}.
\end{equation}
The $n=2$ case is sometime called the ``isodiametric'' inequality.

Therefore, we see that the condition $\area f(\D)=\pi\ (=\area(\D))$ is more
general than $\capa f(\D)= 1\ (=\capa (\D))$, which in turn is more
general than 
$d_n(f(\D))= n^{1/(n-1)}\ (=d_n(\D))$. 

The first main result of this paper is 
the following generalization of the Landau-Toeplitz Theorem.
We consider the following ratios:
\[
\phi_{\ndiam}(r)\defeq \frac{d_{n} (f (r\D))}{d_{n} 
(\D)r}
\qquad\mbox{ and }\qquad
\phi_{\capa}(r)\defeq\frac{\capa(f(r\D))}{\capa
  (r\D)}=\frac{d_\infty(f(r\D))}{d_\infty(r\D)}.
\]
\begin{theorem} \label{thm:diameterssl}
Suppose $f$ is analytic on $\mathbb{D}$. The functions
$\phi_{\capa}(r)$ and $\phi_{\ndiam}(r)$ are increasing and
log-convex. Moreover, they are strictly increasing for $0<r<1$ except
in the special case that $f$ is linear. 
\end{theorem}
It can be checked from the power series expansion of $f$ that the
following limits hold:
\begin{equation}\label{eq:asympt}
\lim_{r\downarrow 0}\phi_{\rad}(r)=
\lim_{r\downarrow 0}\phi_{\ndiam}(r)=\lim_{r\downarrow
  0}\phi_{\capa}(r)=|f^{\prime}(0)|. 
\end{equation}
Hence,
we leave as an exercise to show that Theorem \ref{thm:diameterssl}
implies the following corollary.  
\begin{corollary}\label{cor:diameterssl}
Suppose $f$ is analytic on $\D$ and
$d_n(f(\D))=d_n (\D)$ (or $\capa f(\D)=\capa \D$).   
Then  
\begin{align}
\hbox{\,}&\label{eq:ndiamineq} d_n(f(r\D))
\leq d_n(\D)r,\ \qquad\mbox{for every $r\in(0,1)$}\\  
&\label{eq:capalt} \mbox{(resp.~$\capa f(r\D)\leq (\capa~\D)r$) and}\\
&\label{eq:translt} |f^\prime(0)|\leq 1.
\end{align}
Moreover, equality holds in {\rm
(\ref{eq:ndiamineq})} (resp.~in {\rm (\ref{eq:capalt})}) for some $0<r<1$, or in {\rm (\ref{eq:translt})}, if
and only if 
$f (z)$ is an
Euclidean isometry $a+cz$
for some constants  $a\in \C$, $c\in \T$.
\end{corollary}
\begin{remark}\label{rem:loglog}
It follows from the proof of Theorem \ref{thm:diameterssl}, see
the proof of Lemma \ref{T:DSH}, that $\phi_{\ndiam}(r)$ and
$\phi_{\capa}(r)$ have the stronger property that their logarithms are 
convex functions of $\log r$. This is also how Hadamard's Theorem is
usually phrased.
\end{remark}

\subsection{An Area Schwarz Lemma}\label{ssec:areasl}
As mentioned above the condition that $\area f(\D)=\pi$ is
weaker than $\diam f (\D)=  2$. 
We can prove the following analog of Schwarz's Lemma.

\begin{theorem}[Area Schwarz's Lemma]\label{thm:asl}
Suppose $f$ is analytic on the unit disk $\D$. 
Then the function $\phi_{\area} (r)\defeq  (\pi r^2)^{-1}\area f (r\D)$ is 
strictly increasing for $0<r<1$, except when $f$ is linear, in which
case $\phi_{\area}$ is constant.
\end{theorem}
Moreover, by the power series expansion, $\lim_{r\downarrow
  0}\phi_{\area} (r) =|f^\prime(0)|^2.$ So the following corollary ensues.
\begin{corollary}\label{cor:asl}
Suppose $f$ is analytic on the unit disk $\D$ with $\area f (\D)=\pi$. Then
\begin{align}
&\label{eq:agrowth}
\area f (r\D)\leq \pi r^2 \qquad  \mbox{for every $0<r<1$, and}\\
& \label{eq:ader} |f^\prime(0)|\leq 1. 
\end{align}
Moreover, equality holds in 
{\rm (\ref{eq:agrowth})} for some $0<r<1$, or in {\rm (\ref{eq:ader})},
if and only if $f (z)$ is an Euclidean isometry $a+cz$
for some constants  $a\in \C$, $c\in \T$.
\end{corollary}

One might ask whether $\phi_{\area}(r)$ is also
log-convex as with the 
growth functions $\phi_{\rad}$, $\phi_{\ndiam}$, and $\phi_{\capa}$. 
This is true for univalent functions, but fails in general.
In Section \ref{sec:example} we give an explicit example for which
$\phi_{\area}(r) $ is not log-convex.

\subsection{Structure of the paper and other results}

The structure of the paper is as follows. In Section \ref{sec:trans}
we prove Theorem \ref{thm:diameterssl} about $n$-diameter and
capacity generalizations of Schwarz's Lemma. In Section
\ref{sec:original} we dust off the original approach of Landau and
Toeplitz and show that it can be made to work for $n$-diameters,
hence giving an alternative, more direct proof of Theorem
\ref{thm:diameterssl} for $n$-diameter. In Section \ref{sec:areasl} we
explore an even further generalization of Schwarz's Lemma using area
instead and prove Theorem~\ref{thm:asl}. In Section \ref{sec:example},
however,  
we give an example where log-convexity fails. In Section \ref{sec:perim}, we
formulate a
generalization using perimeter. This is our weakest result because
log-convexity is missing and extra conditions must be imposed on the
range. In Section \ref{sec:hypgeom}, we give some applications of
these Schwarz lemmas to hyperbolic geometry. In particular, we obtain the
global lower bound (\ref{eq:lowerbd}) for the Poincar\'e density on
arbitrary domains. In Section
\ref{sec:modgrowth}, we study bounds on the growth of $|f(z)|$ under
conditions on the image $f(\D)$ that involve diameter instead of
radius.
In Section \ref{sec:highder} we describe a related result of Poukka, obtained
around the same time as the Landau-Toeplitz paper, which involves
higher derivatives. Finally in Section \ref{sec:furtherpb} we state
some open problems.

\section{Higher and transfinite diameter generalizations of
  Schwarz's Lemma} 
\label{sec:trans}

In this section we prove Theorem \ref{thm:diameterssl}. 
The stronger notion of log-convexity turns out to be essential to
prove the sharp result.

\begin{lemma} \label{T:DSH}
For $f$ analytic on $\mathbb{D}$ and $n=2,3,...$, both 
$\phi_{\ndiam}(r)$ and $\phi_{\capa}(r)$ 
 are increasing convex functions of $\log r$, $0<r<1$. 
\end{lemma}

\begin{proof}
Let $f$ be analytic on $\mathbb{D}$; we may assume that $f(0)=0$ and
that $f$ is not linear. It
suffices to prove that $\phi_{\ndiam}(r)$ is an increasing convex
function of $\log r$, since the corresponding result for
$\phi_{\capa}(r)$ then follows by a limit argument. 

So fix $n=2,3,...$, and consider the auxiliary function
\begin{equation}\label{eq:aux}
F_{w_{1},\dotsc ,w_{n}} (z)\defeq d_{n}
(\D)^{-\frac{n(n-1)}{2}}\prod_{j<k} (f (w_{k}z)-f (w_{j}z)),  
\end{equation}
for fixed distinct $w_{1},\dotsc ,w_{n}\in \overline{\D}$.
Then $F_{w_1,...,w_n} (z)=z^{\frac{n (n-1)}{2}}g (z)$, where $g$ is
analytic in $\D$. 
So 
\[
\log\left(r^{-\frac{n
    (n-1)}{2}}\rad F_{w_{1},\dotsc ,w_{n}} 
(r\D)\right)=\max_{|z|<r} \log|g(z)|
\]
is strictly increasing for $0<r<1$, except in the special case when
$g(z)\equiv g(0)$; in fact, by Hadamard's three-circles Theorem it is
also log-convex. Moreover,
for fixed $r\in (0,1)$ we have
\begin{equation}\label{eq:maxrad}
\max_{w_{1},\dotsc ,w_{n}\in \overline{\D}}\rad  F_{w_{1},\dotsc,w_{n}} (r\D)=
\left(\frac{d_{n} (f (r\D))}{d_{n} (\D)} \right)^{\frac{n (n-1)}{2}}.
\end{equation}
So the function 
\begin{equation}\label{eq:phind}
\log\phi_{\ndiam}(r)=\max_{w_{1},\dotsc ,w_{n}\in
  \overline{\D}}\log\left(r^{-\frac{n(n-1)}{2}}\rad  F_{w_{1},\dotsc,w_{n}}
(r\D)\right)^{\frac{2}{n(n-1)}}
\end{equation}  
is the pointwise maximum
of a family of increasing log-convex functions, hence it is increasing
and log-convex for $0<r<1$. This 
implies that $\phi_{\ndiam}(r)$ itself must be increasing and log-convex.
So Lemma \ref{T:DSH} is proved.
\end{proof}

Finally, we will need the following elementary lemma.
\begin{lemma}\label{lem:elem}
Let $f$ be analytic in $\D$ and not linear. Then there is $0<r_0<1$ such that
for $0<r<r_0$,
\[\phi_{\area}(r)>|f^\prime(0)|^2.
 \]
\end{lemma}
\begin{proof}
The statement is clear if $f^\prime(0)=0$. So assume $f^\prime(0)\neq
0$. Then $f$ is one-to-one near the origin and for $r>0$ small
\[
\area f(r\D)=\int_{r\D}|f^\prime(z)|^2d\cH^2(z)=\pi\sum_{n=0}^\infty
n|a_n|^{2}r^{2n}
\]
(where $\cH^2$ is two-dimensional Lebesgue measure).
So $\phi_{\area}(r)=\sum_{n=1}^\infty n|a_n|^{2}r^{2(n-1)}$ is
strictly increasing unless $f$ is linear.
\end{proof}

\begin{proof}[Proof of Theorem \ref{thm:diameterssl}]
We do the proof for $\phi_{\capa}(r)$, since the one for
$\phi_{\ndiam}(r)$ is the same except for the obvious changes, e.g.,
use (\ref{eq:areadn}) in place of (\ref{eq:polya}) below.

By Lemma~\ref{T:DSH} the function $\phi_{\capa}(r)$ is
an increasing convex function of $\log r$. Suppose it fails to be 
strictly increasing. Then by monotonicity it must be constant on an interval
$[s,t]$ for some $0<s<t<1$.
By log-convexity, it then would have to be
constant and equal to $|f^\prime(0)|$ on all of the interval $(0,t)$. 
But, for $0<r<\min\{r_0,t\}$, with $r_0$ as in Lemma~\ref{lem:elem},
\begin{equation*}
|f'(0)|^2 \leq \phi_{\area}(r) = \frac{\area\,f(r\mathbb{D})}{\pi r^2} \leq \phi_{\capa}^2(r) = |f'(0)|^2.
\end{equation*}
where P\'olya's inequality (\ref{eq:polya}) has been used. Therefore,
$\phi_{\area}(r)$ is constant on $(0,t)$. So, by Lemma~\ref{lem:elem}, $f$
must be linear.
\end{proof}

\section{The original Landau-Toeplitz approach}
\label{sec:original}

In this section we revive the original method of Landau and Toeplitz.
We show that it can be used to give a direct proof of
Theorem \ref{thm:diameterssl} for $n$-diameters. However, it seems that
for capacity one really needs to use log-convexity and P\'olya's
inequality.

The proof hinges on the following lemma.
\begin{lemma}\label{lem:fp}
Suppose $g$ is analytic in $\D$, $0<r<1$, $|w|=r$, 
\[
w=g(w)\mbox{ and }r=\max_{|z|=r}|g(z)|.
\]
Then, $\ima g^\prime(w)=0$.
\end{lemma}
\begin{proof}
Actually, the stronger conclusion $g^\prime(w)\geq 0$ is geometrically
obvious because when $g^\prime(w)\neq 0$, the map $g$ is very close to
the rotation-dilation centered at $w$ given by $\zeta\mapsto
w+g^\prime(w)(\zeta-w)$. Since $g$ can't rotate 
points inside $D(0,|w|)$ to  points outside, the derivative must be
positive.

For the sake of rigor, we instead give a ``calculus'' proof of the weaker
statement, along the lines of the original paper
of Landau and Toeplitz, which they credit
to F. Hartogs. 

For $\theta\in \R$ introduce
\[
\phi(\theta)\defeq
|g(we^{i\theta})|^2=g(we^{i\theta})\overline{g(we^{i\theta})}.
\]
The function $g^\star(z)\defeq\overline{g(\bar{z})}$ is also analytic
in $\D$, and $\phi$ may be written
\[
\phi(\theta)=g(we^{i\theta})g^\star(\bar{w}e^{-i\theta}),
\]
enabling us to compute $\phi^\prime(\theta)$ via the product and chain
rules. We get routinely,
\[
\phi^\prime(\theta)=-2\ima \left[we^{i\theta}g^\prime(we^{i\theta})\overline{g(we^{i\theta})}\right]
\]
and setting $\theta=0$,
\[
\phi^\prime(0)=-2\ima
\left[wg^\prime(w)\overline{g(w)}\right]
=-2\ima \left[wg^\prime(w)\overline{w}\right]
=-2|w|^2\ima g^\prime(w).
\]
Since $\phi$ realizes its maximum over $\R$ at $\theta=0$, we have
$\phi^\prime(0)=0$, so the preceding equality proves Lemma \ref{lem:fp}. 
\end{proof}

The following fact will also be important in the sequel.
\begin{fact}\label{fact}
Given $n$ points $\{w_j\}_{j=1}^n\subset
\overline{\D}$, 
\[
\prod_{j<k}|w_j-w_k|\leq n^{\frac{n}{2}}
\]
with equality if and only if, after relabeling, $w_j=u\alpha^j$ for
some $u\in \T$, where $\alpha^{j}$
are the $n$-th roots of unity: i.e.,
$\alpha^j\defeq \exp(i(2\pi j)/n)$.
\end{fact}

We briefly sketch here why this is so. 
Recall that given $n$ complex numbers
$\{w_{j} \}_{j=1}^{n}$, one may form the Vandermonde matrix
$V_n\defeq [w_{j}^{k-1}]_{j,k=1}^{n}$, and that 
\begin{equation}\label{eq:detVn}
\det V_n = \prod_{1\leq j<k\leq n} (w_{k}-w_{j}).
\end{equation}
Indeed, $\det V_n$ is a polynomial of degree at most $n-1$ in $w_n$,
vanishing at $w_1,\dots,w_{n-1}$ with coefficient of $w_n^{n-1}$ equal
to $\det V_{n-1}$, so that (\ref{eq:detVn}) follows by induction.
Hadamard's inequality states that for every $n\times n$ matrix
$A=[a_{jk}]$ with complex entries:
\[
\displaystyle \vert\det(A)\vert\leq
\prod_{j=1}^n\left(\sum_{k=1}^n\vert
a_{jk}\vert^2\right)^\frac{1}{2},
\]
with equality if and only if the rows of $A$
are orthogonal.

Since all the entries of the matrix $V_n$ are bounded by $1$ in modulus,
we find that 
\begin{equation}\label{eq:adetV}
\vert\det(V_n)\vert\leq n^{\frac{n}{2}} 
\end{equation}
with equality if and only if the rows of $V_n$ are orthogonal, i.e., if
and only if 
\[
0=\sum_{k=1}^{n}
w_{j}^{k-1}\overline{w_{l}}^{k-1}=\frac{(w_{j}\overline{w_{l}})^{n}-1}{w_{j}\overline{w_{l}}-1}  
\]
whenever $j\neq l$. Fact \ref{fact} then follows.

It follows that 
the $n$-diameter of $\D$ is $d_{n}(\D)=n^{1/ (n-1)}$ (which strictly
decreases to 
$d_{\infty} (\D)=1$), and it is attained exactly at the $n$-th roots
of unity (modulo rotations). 
\bigskip

\begin{proof}[Landau-Toeplitz-type proof of Theorem \ref{thm:diameterssl} for $n$-diameter]
Consider as before the auxiliary function
$F_{w_{1},\dotsc ,w_{n}} (z)$ defined in (\ref{eq:aux}),
for fixed distinct $w_{1},\dotsc ,w_{n}\in \overline{\D}$.
Then $F_{w_1,...,w_n} (z)=z^{\frac{n (n-1)}{2}}g (z)$, where $g$ is analytic in $\D$
and 
\begin{equation}\label{eq:gn}
|g (0)|=\left|\frac{f^{\prime} (0)}{d_{n} (\D)} \right|^{\frac{n (n-1)}{2}}
\prod_{j<k}|w_{k}-w_{j}|.
\end{equation}
As shown above, in the proof of Lemma \ref{T:DSH}, the function 
$r\mapsto (d_{n}(\D)r)^{-1}d_{n} (f (r\D))$ is increasing for $0<r<1$.
Assume that it is not strictly increasing. 
Then we can find $0<s<t<1$ so that it constant on $[s,t]$. By
(\ref{eq:maxrad}) we can find 
distinct $w_{1}, \dotsc , w_{n}\in \overline{\D}$ so that
\begin{equation}\label{eq:attained}
\left(\frac{d_{n} (f (s\D))}{d_{n}(\D)s}\right)^{\frac{n
(n-1)}{2}}= s^{-\frac{n 
(n-1)}{2}}\rad F_{w_{1},\dotsc ,w_{n}} (s\D). 
\end{equation}
But (\ref{eq:maxrad}) also implies that
\[
\left(\frac{d_{n} (f (r\D))}{d_{n}(\D)r}\right)^{\frac{n
(n-1)}{2}}\geq r^{-\frac{n 
(n-1)}{2}}\rad F_{w_{1},\dotsc ,w_{n}} (r\D) 
\]
for every $0<r\leq t$. In particular, letting $r=t$ and
by Schwarz's Lemma (Theorem \ref{thm:sl}) applied to $g$ (for this
choice of $w_j$'s), we find that
$r^{-\frac{n(n-1)}{2}}\rad F_{w_{1},\dotsc ,w_{n}} (r\D)$ is constant 
for $0<r\leq t$; hence by (\ref{eq:attained}) and the monotonicity of
$\phi_{\ndiam}(r)$,
\begin{equation}\label{eq:constant}
\left(\frac{d_{n} (f (r\D))}{d_{n}(\D)r}\right)^{\frac{n
(n-1)}{2}}\equiv r^{-\frac{n 
(n-1)}{2}}\rad F_{w_{1},\dotsc ,w_{n}} (r\D) = \left(\frac{|f^{\prime}
(0)|}{d_{n} (\D)} \right)^{\frac{n (n-1)}{2}} 
\prod_{j<k}|w_{j}-w_{k}|,
\end{equation}
for $0<r<t$.
In particular, either $f$ is constant or $f^{\prime} (0)\neq
0$. In what follows assume 
$f$ is not constant.

We have 
\[
|f^{\prime}
(0)|^{\frac{n (n-1)}{2}}=\lim_{z\rightarrow 0}\frac{1}{(d_{n} (\D))^{\frac{n (n-1)}{2}}}\prod_{j<k}\left|
\frac{f (\alpha^{k}z)-f (\alpha^{j}z)}{z}\right|\leq
\lim_{r\rightarrow 0}\left(\frac{d_{n} (f (r\D))}{rd_{n} (\D)}
\right)^{\frac{n (n-1)}{2}}, 
\]
and so from (\ref{eq:constant}),
\[
\prod_{j<k}|w_{j}-w_{k}|\geq (d_{n} (\D))^{\frac{n (n-1)}{2}},
\] 
which implies that $w_{j}=u\alpha^{j}$ for some $u\in \T$ by Fact \ref{fact}.
By a rotation, we may take $u=1$.

Therefore, we find that,
for all $z\in \D$
\begin{equation}\label{eq:alphas}
F_{\alpha^1,...,\alpha^n}(z)= d_{n}
(\D)^{-\frac{n(n-1)}{2}}\prod_{j<k} (f (\alpha^{k}z)-f
(\alpha^{j}z))=c(zf^\prime(0))^{\frac{n(n-1)}{2}},
\end{equation}
where $c$ is a constant with $|c|=1$.
In particular, notice that $f (z\alpha^{k})-f (z\alpha^{j})=0$ if and
only if $z=0$.

Now, fix $0<|z|=r<t$ and consider the function 
\[
h_{z} (\zeta)\defeq \prod_{k=1}^{n-1}\frac{f (\zeta z )-f
  (z \alpha^{k})}{(1-\alpha^{k}) zf^{\prime}
  (0)}\prod_{1<j<l\leq n-1}\frac{f (z \alpha^{l})-f (z \alpha^{j})}{(\alpha^{l}-\alpha^{j}) zf^{\prime} (0)},
\]
which is analytic for $\zeta \in \overline{\D}$. 
Then by (\ref{eq:alphas}) and Fact \ref{fact},
\[
|h_{z} (1)|=1 \geq  \sup_{|\zeta |<1}|h_{z} (\zeta )|.
\]
Note that, 
\[
h_z^\prime(\zeta)=h_z(\zeta)z f^\prime(\zeta z )\sum_{k=1}^{n-1}\frac{
1}{f(z \zeta)-f(z \alpha^k)}. 
\]
By Lemma \ref{lem:fp} applied to $h_z (\cdot) /h_z(1)$, and the Open-Mapping
Theorem there is a real constant 
$A$ so that
\begin{equation}\label{eq:fprime}
z  f^\prime(z)\sum_{k=1}^{n-1}\frac{1}{f(z)-f(z \alpha^k)}=A,
\end{equation}
for $z\in t\D\setminus\{0\}$.

To show that (\ref{eq:fprime}) implies $f$ is linear, we may suppose
$f(0)=0$, $f'(0)=1$. In the $n=2$ case (the one considered in the
Landau-Toeplitz paper), the end-game is much simpler. Here, in the
general case, we proceed as follows.
If $f$ is not linear, we may
write $f(z)=z+a_p z^p +\dots$ where $a_p \ne 0$, and $p\geq 2$. Then
\[
\frac{1}{f'(z)}=1 - p a_p z^{p-1}+\dots,
\]
and
\[
\sum_{k=1}^{n-1} \frac{z}{f(z)-f(z\alpha^k)}=
\sum_{k=1}^{n-1} \frac{1}{1-\alpha^k} - \sum_{k=1}^{n-1} 
\frac{1-\alpha^{kp}}{(1-\alpha^k)^2}a_p z^{p-1}+\dots.
\]
This and (\ref{eq:fprime}) imply that 
\[
A=\sum_{k=1}^{n-1}
\frac{1}{1-\alpha^k} ~~~\hbox{ and } ~~~ pA=\sum_{k=1}^{n-1}
\frac{1-\alpha^{kp}}{(1-\alpha^k)^2}.
\]
Recall that $A$ is real and 
$\rea (1/(1-\alpha^k)) =\frac{1}{2}$, so that $A=(n-1)/2.$

For $1\le j \le n$, by Fubini,
\begin{equation}\label{eq:lemma}
\sum_{k=1}^{n-1}\frac{1-\alpha^{jk}}{1-\alpha^k}
=\sum_{k=1}^{n-1} \sum_{q=0}^{j-1}
\alpha^{qk}=n-1+\sum_{q=1}^{j-1}\left[\frac{1-\alpha^{qn}}{1-\alpha^q}
-1 \right] =n-j
\end{equation}
since $\alpha^n=1$.
So, if $1\le p \le n$, using Fubini, (\ref{eq:lemma}) and the
definition of $A$, we get
\begin{equation}\label{eq:pa}
\sum_{k=1}^{n-1}\frac{1-\alpha^{kp}}{(1-\alpha^k)^2}=
\sum_{k=1}^{n-1}\sum_{j=0}^{p-1}\frac{\alpha^{jk}}{1-\alpha^k}=A+\sum_{j=1}^{p-1}[A-(n-j)]=pA-(n-\frac{p}{2})(p-1).
\end{equation}
The earlier identification of $pA$ then leads to the conclusion that
\begin{equation}\label{eq:contra}
pA=pA-(n-\frac{p}{2})(p-1),
\end{equation}
which is a contradiction for $2\le p \le n$. If $p>n$ and $p\equiv p'$ mod $n$,
with $p'\le n$, then $\alpha^{kp}=\alpha^{kp'}$ so (\ref{eq:pa}) again
shows that
\[
pA=p'A-(n-\frac{p'}{2})(p'-1)
\]
which is impossible. Thus, the assumption that $f$ is not linear is untenable.
\end{proof}

\section{Area generalization of Schwarz's Lemma}
\label{sec:areasl}

In this section we prove 
Theorem \ref{thm:asl}. This requires some preliminaries.

Let $f$ be analytic and non-constant in a neighborhood of $\overline{\D}$. 
For every $w\in f(\D)$ let 
\[
Z(w)\defeq
\{z_j(w)\}_{j=1}^{N(w)}
\] 
be the set of points in
$f|_{\D}^{-1}(w)$ of minimum modulus. Note that $0<N(w)<\infty$.
\begin{claim}
The function $w\mapsto N(w)$ is Borel measurable on $f(\D)$.
\end{claim}
\begin{proof}
Let $C=\{z\in \D: f^\prime(z)=0\}$ be the set of critical points,
which is finite,
and let $P=f(C)$ be the finite post-critical set. Pick $w_0\in
f(\D)\setminus P$. It is enough to
check Borel measurability of $N$ near $w_0$. By the
argument principle, there is a small disk $D$ centered at $w_0$ and
there are $M$ branches of the inverse of $f$ such
that $f|_{\D}^{-1}(w)=\{\zeta_1(w),\dots,\zeta_M(w)\}$ for every $w\in
D$; see Theorem on p.~238 of \cite{gamelin2001}. Now, upon relabeling,
$Z(w_0)=\{\zeta_j(w_0)\}_{j=1}^{M^\prime}$ for some $M^\prime\leq
M$. Moreover, restricting to a smaller disk $D^\prime\subset D$
centered at $w_0$, we can assume using continuity of the branches that
$Z(w)\subset \{\zeta_j(w)\}_{j=1}^{M^\prime}$ for every $w\in
D^\prime$. Since each $\zeta_j$ is analytic, standard results show 
that $N$ is Borel measurable.
\end{proof}

Now consider the set 
\[
E\defeq \cup_{w\in f(\D)}Z(w).
\]
\begin{claim}
The set $\D\setminus E$ is open.
\end{claim}
\begin{proof}
Pick $z_0\in \D\setminus E$. Let $w_0=f(z_0)$. Then we can find
$z_0^\star\in f^{-1}(w_0)$ such that
\[
|z_0^\star|<|z_0|.
\]
Assume that $z_0$ is of order $m-1$ and $z_0^\star$ of order
$m^\star-1$. Then, by the argument principle, there are disks
$D=D(z_0,\epsilon)$ and $D^\star=D(z_0^\star,\epsilon)$ with
radius $0<\epsilon<(|z_0|-|z_0^\star|)/3$ small enough so that for every
point $z\in D$, the value $w=f(z)$ is close enough to $w_0$ to have at
least one preimage $z^\star$ in $D^\star$. Thus $D\subset \D\setminus E$.
\end{proof}

Below we will need the following ``Non-Univalent Change
of Variables Formula''. 
\begin{theorem}[Theorem 2, p.~99 of
    \cite{evans-gariepy1992}]\label{thm:chgevar} 
Let $f:\R^n\rightarrow \R^m$ be Lipschitz, $n\leq m$. Then for each
integrable $g:\R^n\rightarrow \R$,
\[
\int_{\R^n}g (x)Jf (x)dx=\int_{\R^m}\left[\sum_{x\in f^{-1} (y)}g(x)
\right]d\cH^{n} (y),
\]
whenever either side converges, and where $Jf$ is the Jacobian
(determinant) of $f$. 
\end{theorem}

\begin{proof}[Proof of Theorem \ref{thm:asl}]
Fix $0<r<1$ and consider the integral
\begin{equation}\label{eq:areaint}
\int_\D\frac{\chi_{r\D\cap E}(z)}{N(f(z))}|f^\prime(z)|^2d\cH^2(z).
\end{equation}
By the non-univalent change of variables (Theorem \ref{thm:chgevar}) 
applied to the function
\[g(z)=\frac{\chi_{r\D\cap E}(z)}{N(f(z))}:\D\rightarrow \R,\] it equals
\begin{eqnarray*}
\int_{f(\D)}\frac{1}{N(w)}\sum_{z\in f^{-1}(w)}\chi_{r\D\cap E}(z)d\cH^2(w) & = & 
\int_{f(r\D\cap E)}\frac{1}{N(w)}\sum_{z\in f^{-1}(w)}\chi_{r\D\cap Z(w)}(z)d\cH^2(w)\\
& = &  \int_{f(r\D\cap E)}\frac{1}{N(w)}\sum_{z\in Z(w)}1d\cH^2(w)
  \\
& = & \area f(r\D\cap E)\\
& = & \area f(r\D).
\end{eqnarray*}
The last equality holds because $w\in f(r\D)$ if
and only if
$Z(w)\subset r\D$ if and only if $w\in f(r\D\cap E)$.

Thus,
\[
\area f(r\D)=\int_\D\frac{\chi_{r\D\cap
E}(z)}{N(f(z))}|f^\prime(z)|^2d\cH^2(z), 
\]
for $0<r<1$. In particular, the function $A (r)\defeq \area f (r\D)$
is absolutely continuous. 

By Fubini,
\[
A (r)=\int_\D\frac{\chi_{r\D\cap
E}(z)}{N(f(z))}|f^\prime(z)|^2d\cH^2(z)=
\int_0^{r}\int_{s\T}\frac{\chi_{
    E}(z)}{N(f(z))}|f^\prime(z)|^2|dz|ds.
\]
So we have
\[
\frac{dA
(r)}{dr}=\int_{r\T}\frac{\chi_{E}(z)}{N(f(z))}|f^\prime(z)|^2|dz|. 
\]
By the
Cauchy-Schwarz inequality, 
\begin{equation*}
\frac{dA(r)}{dr} \geq 
\frac{1}{2\pi r}\left(\int_{r\T}\frac{\chi_{
    E}(z)}{N(f(z))}|f^\prime(z)||dz| \right)^2. 
\end{equation*}
Again, by the non-univalent change of variables (Theorem
\ref{thm:chgevar}) applied to the 
function
\[g(z)=\frac{\chi_{E}(z)}{N(f(z))}:r\T\rightarrow\R,\] we obtain
\begin{eqnarray*}
\int_{r\T}\frac{\chi_{
    E}(z)}{N(f(z))}|f^\prime(z)||dz|
    & = & \int_{f(r\T\cap E)}\frac{1}{N(w)}\#(
    r\T\cap Z(w))d\cH^1(w)\\
& = & \length f(r\T\cap E)=\length \bd f(r\D).
\end{eqnarray*}
The last equality holds because $w\in \bd f (r\D)$ if and only if
$w\in f (r\T\cap E)$. 
Thus, writing $L (r)\defeq \length \bd f (r\D)$, we have shown that
\[
\frac{dA(r)}{dr}\geq \frac{L (r)^{2}}{2\pi r}.
\]
The isoperimetric inequality \cite{lax1995} says that, for planar domains,
\[
4\pi\area
\Om\leq (\length \bd\Om)^{2}.
\]
So, we have
\begin{equation}\label{eq:au}
\frac{dA(r)}{dr}\geq \frac{2A(r)}{r}.
\end{equation}
Now consider the function $\phi_{\area} (r)$ defined in the statement
of Theorem \ref{thm:asl}. We have
shown that it is absolutely continuous and its derivative is
\[
\frac{d\phi_{\area} (r)}{dr}=-2\pi^{-1}r^{-3}A (r)+ (\pi
r^2)^{-1}\frac{dA (r)}{dr}= (\pi r^2)^{-1}\left(\frac{dA (r)}{dr}-\frac{2A
(r)}{r} \right)\geq 0 
\]
by (\ref{eq:au}).
Therefore, $\phi_{\area} (r)$ is an increasing function of $r$.

If $\phi_{\area} (r)$ is not strictly increasing, then there is
$0<s<t<1$ such that $\phi_{\area} (r)=c$ for every $s\leq r\leq
t$. This implies that $\phi_{\area}^{\prime}(r)\equiv 0$ on
$[s,t]$. Hence, 
\[
\frac{dA(r)}{dr}\equiv \frac{L (r)^{2}}{2\pi r}\equiv\frac{2A(r)}{r}
\]
on $[s,t]$. So the extremal case in the isoperimetric inequality shows
that $f (r\D)$ is a disk for $s\leq r\leq t$, with area $\pi r^{2}$.
Hence $\rad f (r\D)\equiv r$ on $[s,t]$, and by Theorem \ref{thm:sl},
we conclude 
that $f$ must be linear. 
\end{proof}
We leave the proof of Corollary \ref{cor:asl} as an exercise for the reader.

\section{A counter-example to log-convexity}\label{sec:example}

Notice that $\log\phi_{\area}(r)$ is a convex function of $\log r$ if
and only if $\log A(r)$ is. Also,
$\log A(r)$ is log-convex for all univalent functions. 
In fact, write $f(z)=\sum_{n=0}^\infty a_n z^n$. If $f$ is univalent,
$A(r)=\sum_{n=0}^\infty n|a_n|^2 r^{2n}$. Then, by straight differentiation, 
$ \log \sum_{n=1}^\infty n|a_n|^2 e^{2nx}$ has non-negative second 
derivative if and only if
\[
\sum_{n,k=1}^\infty 4(n k^3-n^2 k^2)|a_n|^2 |a_k|^2  e^{(2n+2k)x} \geq 0.
\]
If we switch $n$ and $k$ and add the results, it doesn't affect the 
truth of non-negativity, so the above will be non-negative if
\[
n k^3 +k n^3 \geq 2n^2 k^2
\]
dividing by $n^2 k^2$ it suffices that
$k/n+n/k \ge 2$, which is true. In fact equality occurs if and only if 
$n=k$, and hence we have a strictly positive second derivative unless
$f(z)=cz^m$, and by univalence, unless $f$ is linear.

However, as the following example shows, neither $\log 
A(r)$ nor even $A(r)$ is log-convex in general.
\begin{example}
We study the function
\[f(z)=\exp\left(ic\log\left(\frac{1+z}{1-z}\right)\right)\]
with $c>0$, which is a universal cover of $\D$ onto the annulus
$\{e^{-\pi c/2}<|z|<e^{\pi c/2}\}$. To compute $A(r)$ we first apply the
conformal map $\psi(z)=i\log\left(\frac{1+z}{1-z}\right)=u(z)+iv(z)$ which
sends $r\D$ into an oval contained in the vertical strip
$\{|u|<\pi/2\}$. We then notice that $f(r\D)\setminus
(-\infty,0)$ is covered by the restriction of $e^{\pi z}$
to the part of the oval which is in $\{|v|<\pi/c\}$.
So a computation shows that $f$ is univalent on $r\D$ for
$r<\tanh(\pi/(2c))$ and that for $\tanh(\pi/(2c))\leq r<1$ 
\[
A(r)=\int_0^\pi2\sinh\left(2c\arccos\left(\frac{1-r^2}{1+r^2}\cosh(t/c)\right)\right)dt. 
\]
Writing $A_c(r)$ for $A(r)$ to emphasize the dependence on the
parameter $c$, we then study the asymptotics as $c\downarrow 0$. We find
that for $x\in (0,1)$,
\[
\lim_{c\downarrow
  0}\frac{A_c(e^{-x\log\coth(\pi/(2c))})-2\pi\sinh(c\pi)}{4c^2}=-\int_0^x\frac{\arcsin
  u}{u}du. 
\]
But the right hand-side is a strictly concave function of $x\in(0,1)$,
since its derivative is {\em minus} the strictly increasing
function $x^{-1}\arcsin x$. Thus, for $c>0$ sufficiently small, 
$A_c(r)$ cannot be log-convex.
\end{example}

\section{Perimeter generalization of Schwarz's Lemma}\label{sec:perim}

The results are not as strong when considering the notion of perimeter. 
\begin{proposition}\label{prop:lengthsl}
Suppose $\Om$ is a simply-connected domain. 
If $F$ is a one-to-one analytic map of $\D$ onto $\Om$, then
$r\mapsto r^{-1}\length \bd (F (r\D))$ is strictly increasing, unless
$F$ is linear. Hence, if $\bd\Om$ is a
Jordan curve with Euclidean length at most
$2\pi$, then 
\begin{align}
\hbox{\,}
&\label{eq:lengthineq} \length \bd (F (r\D)) \leq 2\pi r\qquad \mbox{for
every $0<r<1$},\\
&\label{eq:Fderle1} |F^\prime(0)|\leq 1.
\end{align}
Moreover, equality holds in {\rm
(\ref{eq:lengthineq})} for some $0<r<1$, or in {\rm
(\ref{eq:Fderle1})}, if and only if 
$F (z)$ is a Euclidean isometry $a+cz$ for some $a\in \C$, $c\in \T$.
\end{proposition}

\begin{corollary}\label{cor:lengthsl}
Suppose $\Om$ is a simply-connected region in $\C$ and $\bd\Om$ is a
Jordan curve with Euclidean length at most
$2\pi$. If $f$ is analytic on $\D$ with values in $\Om$, then
$|f^{\prime} (0)|\leq 1$ and equality holds if and only if $f (z)$ is
a Euclidean isometry $a+cz$
for some $a\in \C$, $c\in \T$. 
\end{corollary}

\begin{remark}\label{rem:known}
The bound (\ref{eq:Fderle1}) also follows from the well-known result
that the $H^{1}$-norm of
$F^{\prime}$ is the length of the boundary of the image, together with the
mean-value inequality
\[
|F^{\prime}(0)|\le \int |F^{\prime}(re^{it})| \frac{dt}{2\pi}.
\]
The right-hand side converges to the $H^{1}$-norm as $r\uparrow 1$.
Also the right side above increases
with $r$, so that $r$ times the right side, is bounded above by $r$ times
the right side evaluated at $r=1$, and that gives (\ref{eq:lengthineq}).
\end{remark}
\begin{remark}\label{rem:sqtrick}
The same ``square root trick'' used in the proof of Proposition
\ref{prop:lengthsl} below can be
used to prove the isoperimetric inequality; see \cite{carleman1921}
and \cite{duren1983}
exercise 3, page 25. 
\end{remark}
\begin{remark}\label{rem:macgreg}
{\em Lower} bounds for area and perimeter of image disks can be found in a
paper of MacGregor \cite{macgregor1964} and they involve the
derivative at the origin.
\end{remark}

\begin{proof}[Proof of Proposition \ref{prop:lengthsl}]
Let $G (z)\defeq \sum_{n=0}^{\infty}b_{n}z^{n}$ be an analytic square
root in $\D$ of the zero-free function $F^{\prime}$. Then $F^{\prime}
(0)=G^{2} (0)$ and 
\begin{eqnarray}
r^{-1}\length \bd (F (r\D)) & = & r^{-1}\int_{|z|=r}|F^{\prime}
(z)||dz|\notag\\  
& = & r^{-1}\int_{|z|=r}|G (z)|^{2}|dz|\notag\\
& = & 2\pi \sum_{n=0}^{\infty}|b_{n}|^{2}r^{2n},\label{eq:dagger}
\end{eqnarray}
which is strictly increasing for $0<r<1$ unless $b_{n}=0$ for all $n\geq
1$, i.e., 
unless $F$ is linear. The rest follows straightforwardly.

Also, the isoperimetric inequality \cite{lax1995} says that 
\[
4\pi\area
\Om\leq (\length \bd\Om)^{2}.
\]
Therefore, $\area \Om\leq \pi$ and by
Corollary \ref{cor:asl}, $|F^{\prime} (0)|\leq 1$, with equality if
and only if $F(z)=a+cz$ identically for some $a\in \C$, $c\in \T$.
\end{proof}
\begin{proof}[Proof of Corollary \ref{cor:lengthsl}]
If $f$ is analytic on $\D$ with values in $\Om$, let $F$ be the
Riemann map  of $\D$ onto $\Om$ with $F (0)=f (0)$ and $F^{\prime} (0)>0$.
Then $g\defeq F^{-1}\circ f$ is a self-map of the disk which fixes the
origin and $f=F\circ g$. So
\[
|f^{\prime} (0)|=F^{\prime} (0)|g^{\prime} (0)|\leq F^{\prime} (0)\leq 1
\]
with equality if and only if $g$ is a rotation, i.e., if and only if
$f (z)=F (cz)$ for some $c\in \T$. Now apply Proposition
\ref{prop:lengthsl}. 
\end{proof}

\section{Applications to hyperbolic geometry}
\label{sec:hypgeom}

The {\sf hyperbolic metric} on $\mathbb D$ is
$$\rho_{\mathbb D}(z)|dz|\defeq \frac{|dz|}{1-|z|^2}.$$
So $\rho_\D(z)\geq 1$ for every $z\in\D$ with equality when $z=0$.

The associated hyperbolic distance function is
$$h_{\mathbb D}(z,w)\defeq \tanh^{-1}\left|\frac{z-w}{1-\bar w z}\right|.$$
The hyperbolic disk with hyperbolic center $c$ and hyperbolic radius $R>0$ is
$$D_{\mathbb D}(c,R) = \{z:h_{\mathbb D}(z,c)<R\}.$$
The closed hyperbolic disk $\bar D_{\mathbb D}(c,R)$ is defined
similarly. For $0<r<1$ the Euclidean disk $r{\mathbb D}$ is the
hyperbolic disk $D_{\mathbb D}(0,R)$, where $R=\tanh^{-1}r$, or
$r=\tanh R$. 

A region $\Omega$ in $\mathbb{C}$ is {\sf hyperbolic} if
$\mathbb{C}\setminus\Omega$ contains at least two points. If
$\Omega$ is a hyperbolic region and $f:{\mathbb D} \to
\Omega$ is an analytic covering, then the density $\rho_\Omega$ of the
hyperbolic metric $\rho_\Omega(z)|dz|$ on $\Omega$ is defined so that 
\begin{equation} \label{E:DHM}
\rho_\Omega(w)|dw| = \rho_\Omega(f(z))|f'(z)||dz| = \rho_{\mathbb{D}}(z)|dz|.
\end{equation}
This defines the hyperbolic density $\rho_\Omega$ independent of the
covering. Let $h_\Omega$ be the associated hyperbolic distance
function on $\Omega$. Open and closed hyperbolic disks in $\Omega$ are
defined in the standard way. If $f:{\mathbb D}\to\Omega$ is an
analytic covering with $f(0)=c$ and $R>0$, then $f(D_{\mathbb
  D}(0,R))=D_\Omega(c,R)$ with the similar result for closed
hyperbolic disks. 

By Schwarz's Lemma (Theorem \ref{thm:sl}) and the Monodromy Theorem,
the following monotonicity holds 
\[
\tilde{\Omega} \subset \Omega\qquad\Longrightarrow\qquad
\rho_{\tilde{\Omega}}(z) \geq\rho_{\Omega}(z) \qquad\forall
z\in\tilde{\Omega}.  
\]
So for $z\in\Omega$, by choosing $\tilde{\Omega}$ to be the largest
Euclidean disk centered at $z$ contained in $\Omega$, one gets the
following upper bound for hyperbolic density:
\[
\rho_\Omega(z)\leq \frac{1}{\dist(z,\bd\Omega)}.
\]
In \cite[p.~16]{ahlfors1973},
Ahlfors states that it is a 
much harder problem to find lower bounds. Theorem
\ref{thm:arealb} (b) below shows, in particular, that the 
geometric lower bound of 
\begin{equation}\label{eq:lowerbd}
\rho_\Omega(z)\geq \sqrt{\pi/\area(\Omega)}
\end{equation}
holds for every region $\Omega$ and every $z\in\Omega$.

\begin{theorem}\label{thm:arealb}
Suppose $\Omega$ is a hyperbolic region. Then for each $c\in \Omega$
and $R>0$, the function $R \mapsto \left(\pi  \tanh^2(R)\right)^{-1}
{\rm Area}\, D_\Omega(c,R)$ is strictly increasing except when $\Omega$ is
a Euclidean disk with center $c$. If ${\rm Area}\,\Omega\leq\pi$, then
\begin{enumerate}
\item[(a)] for each $c\in\Omega$ and all $R>0$, 
\[\area D_\Omega(c,R)\leq \pi \tanh^2 R\] with equality if and only if $\Omega$
is a Euclidean disk with center $c$ and radius $1$; and 
\item[(b)] for each $c\in\Omega$, \[1 \leq \rho_\Omega(c)\] and equality
holds if and only if $\Omega$ is a Euclidean disk with center $c$ and
radius $1$.
\end{enumerate}
\end{theorem}

\begin{proof} Fix $c\in\Omega$ and let $f:{\mathbb D}\to \Omega$ be an
  analytic covering with $f(0)=c$. Since $f(r\mathbb{D}) =
  D_\Omega(c,R)$, where $r= \tanh R$, 
\begin{equation*}
\frac{{\rm Area}\, D_\Omega(c,R)}{\pi \tanh^2 R} = \frac{{\rm Area}\,
  f(r\mathbb{D})}{\pi r^2}. 
\end{equation*}
Theorem~\ref{thm:asl} implies that this quotient is strictly
increasing unless 
$f$ is linear, or equivalently, $\Omega$ is a disk with center $c$. 

If ${\rm Area}\, \Omega \leq \pi$, then parts (a) and (b) follow from
Corollary~\ref{cor:asl}; note that $\rho_\Omega(c) = 1/|f'(0)|$. 
\end{proof}

Analogous theorems can be formulated for logarithmic capacity and $n$-diameter.

\section{Modulus growth bounds}
\label{sec:modgrowth}

In view of the bound on the growth of the modulus in Schwarz's Lemma, it
is natural to ask whether a similar statement holds in the context of
`diameter'. 
We offer the following result.

\begin{theorem}\label{thm:don}
Suppose $f$ is analytic on the unit disk $\D$ and $\diam f(\D)\leq 2$.
Then for all $z\in \D$
\begin{equation}\label{eq:don}
|f(z)-f(0)| \leq |z|\frac{2}{1+\sqrt{1-|z|^2}}.
\end{equation}
Moreover, equality holds in {\rm (\ref{eq:don})} 
at some point in $\D\setminus\{0\}$  
if and only if $f$ is a linear
fractional transformation of the form
\begin{equation}\label{eq:lft}
f(z)=c\frac{z-b}{1-\overline{b}z}+a
 \end{equation}
for some constants $a\in \C$, $b\in \D\setminus\{0\}$ and $c\in \T$.
\end{theorem}
\begin{remark}
In Schwarz's Lemma, equality in (\ref{eq:growth}) at some point
in $\D\setminus\{0\}$ holds if and only if equality holds at every
point $z\in \D$. This is not true any more in Theorem
\ref{thm:don}. Namely, when
$f$ is the linear fractional transformation in (\ref{eq:lft}), then equality in (\ref{eq:don})
occurs only for $z\defeq 2b/(1+|b|^2)$.  
\end{remark}
\begin{remark}
Since in (\ref{eq:don}) the origin does not play a special role, we can rewrite that inequality
more symmetrically as follows:
\[|f(z)-f(w)| \le \diam f(\D) \frac{\delta}{1+\sqrt{1-\delta^2}}=\diam 
f(\D)\tanh(\rho/2)  \qquad\forall z,w \in \D\]
where $\delta=\delta(z,w)\defeq  \left|\frac{z-w}{1-{\bar w}z}\right|$ is the 
pseudohyperbolic distance between $z$ and $w$ and 
$\rho=\rho(z,w)\defeq  (1/2)\log[(1+\delta)/(1-\delta)]$ is the 
hyperbolic distance between $z$ and $w$.

The preceding inequality can also be rewritten using the well-known 
identity
\[1-\left|\frac{z-w}{1-{\bar w}z}\right|^2 = 
\frac{(1-|z|^2)(1-|w|^2)}{|1-{\bar w}z|^2},\]
as
\[|f(z)-f(w)| \le \diam f(\D) \frac{|z-w|}{|1-{\bar w}z| + 
\sqrt{(1-|z|^2)(1-|w|^2)}}.\]
\end{remark}

\begin{proof}[Proof of Theorem \ref{thm:don}] Fix $d\in\D$ such that $f(d)\neq f(0)$. Set 
$$g= c_1 f\circ T +c_2$$
where $T$ is a linear fractional
transformation of $\D$ onto $\D$ such that
$T(x)=d$, $T(-x)=0$, for some $x>0$  and $c_1$, $c_2$ are constants
chosen so that  $g(x)=x$ and $g(-x)=-x$. By elementary algebra
$$T(z)=\frac{d}{|d|}\frac{z+x}{1+xz}$$
where $x\defeq |d|/(1+\sqrt{1-|d|^2})$, 
$$c_1\defeq \frac{2x}{f(d)-f(0)}\qquad \hbox{ and }\qquad
c_2\defeq -x\frac{f(d)+f(0)}{f(d)-f(0)}.$$
Then 
\begin{equation}\label{eq:diam}
\diam {g(\D)}=|c_1|\diam f(\D)\leq
\frac{4}{|f(d)-f(0)|}\frac{|d|}{(1+\sqrt{1-|d|^2})}.
\end{equation}
We now prove that $\diam g(\D) \ge 2$, with equality if and only
if $g(z)\equiv z$.

Set $h(z)\defeq (g(z)-g(-z))/2.$ Then $h(x)=x$
and $h(-x)=-x$. Note also that $h(0)=0$ so that $h(z)/z$ is analytic
in the disk 
and has value $1$ at $x$ and hence by the maximum principle
$\sup_\D |h(z)|= \sup_\D|h(z)/z|\ge 1$, with equality only if $h(z)=z$
for all $z\in \D$. Since, by definition of $h$, $\diam g(\D)\geq 2\sup_\D
|h|$, we see that $\diam g(\D)\geq 2$ and then (\ref{eq:diam})
gives (\ref{eq:don}) for $z=d$.

If equality holds in (\ref{eq:don}) at some point in
$\D\setminus\{0\}$, then that point is an eligible $d$ for the
preceding discussion, and (\ref{eq:diam}) shows that $\diam g(\D)\leq
2$, while we have already shown that $\diam g(\D)\geq 2$. Thus
$\diam g(\D) =2$. Hence $\sup_{z\in\D}|h(z)|=1$ and therefore $h(z)\equiv z$. 
Since $h$ is
the odd part of $g$,  we have
$g^\prime(0)=h^\prime(0)=1$. Thus, 
by the Landau-Toeplitz Theorem \ref{thm:lt} applied to $g$, we find that  $g(z)\equiv g(0)+z$ and thus
\[
f(z)=\frac{1}{c_1}T^{-1}(z)+f(T(0)).
\]
Moreover, equality at $z=d$ in (\ref{eq:don}) says that $|f(d)-f(0)|=2x$, hence $|c_1|=1$.  Since $T$ is a M\"{o}bius transformation of $\D$, namely of the form
\[
T (z)=\eta\frac{z-\xi}{1-\overline{\xi}z}
\]
for some constants $\xi\in \D$ and $\eta\in \T$,
its inverse is also of this form. Therefore, we conclude that $f$ can be written as in (\ref{eq:lft}).

Finally, if $f$ is given by (\ref{eq:lft}), then $2b/(1+|b|^2)\in \D\setminus\{0\}$, and one checks that equality is attained in (\ref{eq:don}) when $z$ has this value and for no other value in $\D\setminus\{0\}$.
\end{proof}

\section{Higher derivatives}
\label{sec:highder}

We finish with a result, due to Kalle Poukka in 1907, which is to be compared with the usual Cauchy estimates that one gets from the maximum modulus.
Interestingly, Poukka seems to have been the first student of 
Ernst Lindel\"{o}f, who is often credited with having founded the
Finnish school of analysis.  
\begin{theorem}[Poukka \cite{poukka1907}]\label{thm:poukka}
Suppose $f$ is analytic on $\D$. Then for all positive integers $n$ we have
\begin{equation}\label{eq:pk}
\frac{|f^{(n)}(0)|}{n!}\leq \frac{1}{2}\diam f(\D).
\end{equation}
Moreover, equality holds in {\rm (\ref{eq:pk})} for some $n$ if and
only if $f(z)=f(0)+cz^n$ for some constant $c$ of modulus $\diam f(\D)/2$.
\end{theorem}
\begin{proof}[Proof (Poukka):]
Write $c_k\defeq f^{(k)}(0)/k!$, so that $f(z)=\sum_{k=0}^\infty c_k
z^k$, for every $z\in \D$. Fix $n\in \N$. For every $z\in\D$ 
\begin{equation}\label{eq:coef}
h(z)\defeq f(z)-f(ze^{i\pi/n})=\sum_{k=1}^\infty c_k(1-e^{i\pi k/n})z^k.
\end{equation}
Fix $0<r<1$ and notice that, by absolute and uniform convergence,
\begin{equation}\label{eq:parseval}
\sum_{k=1}^\infty|c_k|^2|1-e^{i\pi k/n}|^2 r^{2k}=\int_0^{2\pi} |h(re^{i\theta})|^2\frac{d\theta}{2\pi}\leq (\diam f(\D))^2.
\end{equation}
Therefore
\[
|c_k(1-e^{i\pi k/n})|r^{k}\leq \diam f(\D)
\]
for every $0<r<1$ and every $k\in \N$.
In particular, letting $r$ tend to $1$ and then setting $k=n$, we get $2|c_n|\leq\diam f(\D)$, which is (\ref{eq:pk}). 

If equality holds here, then letting $r$ tend to $1$ in (\ref{eq:parseval}), we get that all coefficients $c_k(1-e^{i\pi k/n})$ in
(\ref{eq:coef}) for $k\neq n$ must be $0$. Hence, $c_k=0$ whenever $k$
is not a multiple of $n$. Thus, $f(z)=g(z^n)$ for some analytic
function $g$ on $\D$. Moreover, $g^\prime(0)=c_n$ and $\diam
g(\D)=\diam f(\D)$. So, by Theorem \ref{thm:lt}, $g(z)=cz$ for some
constant $c$ with $|c|=\diam g(\D)$, and the result follows.

\end{proof}

\section{Further problems}
\label{sec:furtherpb}

Here we discuss a couple of problems that are related to these ``diameter" questions.

The first problem arises 
when trying to estimate the distance of $f$ from its linearization, $f(z)-(f(0)+f^\prime(0)z)$, to give a ``quantitative" version for the `equality' case in Schwarz's Lemma (Theorem \ref{thm:sl}).
This is done via the so-called Schur algorithm. As before, one considers the function
\[
g(z)\defeq \frac{f(z)-f(0)}{z}
\] 
which is analytic in $\D$, satisfies $g(0)=f^\prime(0)$ and which, by
the Maximum Modulus Theorem, has, say, $\sup_\D|g|\leq 1$. Now let $a\defeq f^\prime(0)$ and post-compose $g$
with a M\"{o}bius transformation of $\D$ which sends $a$ to $0$ to
find that
\[
\frac{g(z)-a}{1-\bar{a}g(z)}=zh(z)
\]
for some analytic function $h$ with $\sup_\D|h|\leq 1$.

Inserting the definition of $g$ in terms of $f$ and solving for $f$ shows that
\[
f(z)-f(0)-az=(1-|a|^2)\frac{z^2h(z)}{1+\bar{a}zh(z)}.
\]
Thus, for every $0<r<1$,
\begin{equation}\label{eq:schur}
\max_{|z|<r}|f(z)-f(0)-f^\prime(0)z|\leq (1-|f^\prime(0)|^2)\frac{r^2}{1-|f^\prime(0)|r} 
\end{equation}
and `equality' holds for at least one such $r$ if and only if $h(z)\equiv a/|a|=f^\prime(0)/|f^\prime(0)|$, i.e., if and only if
\[
f(z)=z\frac{a}{|a|}\frac{z+|a|}{1+|a|z}+b
\]
identically, for constants $a\in\overline{\D}$, $b\in \C$.
%

In the context of this paper, when $f$ is analytic in $\D$ and $\diam f(\D)\leq 2$,
by the Landau-Toeplitz Theorem \ref{thm:lt} and a normal-family argument we see that, for every $\epsilon>0$ and every $0<r<1$, there exists $\alpha>0$ such that: $|f^\prime(0)|\geq 1-\alpha$ implies 
\[
|f(z)-(f(0)+f^\prime(0)z)|\leq \epsilon\qquad \forall |z|\leq r.
\]
However, one could ask for an explicit bound as in (\ref{eq:schur}).
\begin{problem}\em
If $f$ is analytic in $\D$ and $\diam f(\D)\leq 2$, 
find an explicit (best?) function $\phi(r)$ for $0\leq r<1$ so that
\[
|f(z)-(f(0)+f^\prime(0)z)|\leq (1-|f^\prime(0)|)\phi(r)\qquad \forall |z|\leq r.
\]
\end{problem}

Another problem can be formulated in view of Section \ref{sec:hypgeom}.
It is known, see the Corollary to Theorem~3 in \cite{minda-wright1982}, that if $\Omega$ is a bounded convex domain, then the minimum  
\begin{equation}\label{eq:min}
\Lambda(\Omega)\defeq \min_{w\in \Omega}\rho_\Omega(w)
\end{equation}
is attained at a unique point $\tau_\Omega$, which we can call the
{\sf hyperbolic center} of $\Omega$. Also let us define the {\sf
  hyperbolic radius} of $\Omega$ to be 
\[
R_h(\Omega)\defeq\sup_{w\in \Omega}|w-\tau_\Omega|.
\]

Now assume that
$\diam \Omega=2$. Then we know, by the corresponding ``diameter''
version of Theorem \ref{thm:arealb}, that
$\Lambda(\Omega)\geq 1$,  with equality if and only if $\Omega$ is a
disk of radius $1$. In particular, if $\Lambda(\Omega)=1$, then
$R_h(\Omega)=1$. 
\begin{problem}\em
Given $m>1$, find or estimate, in terms of $m-1$,
\[
\sup_{\Omega\in \cA_m} R_h(\Omega)
\]
where $\cA_m$ is the family of all convex domains $\Omega$ with $\diam \Omega=2$ and $\Lambda(\Omega)\leq m$.
\end{problem}

More generally, given an analytic function $f$ on $\D$ such that
$\diam f(\D)\leq 2$, define  
\[
M(f)\defeq\min_{w\in\D}\sup_{z\in\D}|f(z)-f(w)|
\]
and let $w_f$ be a point where $M(f)$ is attained. 
\begin{problem}\em
Fix $a<1$. Find or estimate, in terms of $1-a$,
\[
\sup_{f\in \cB_a}M(f)
\]
where $\cB_a$ is the family of all analytic functions $f$ on $\D$ with
$\diam f(\D)\leq 2$ and \[|f^\prime(w_f)|(1-|w_f|^2)\geq a.\]
\end{problem}

Similar questions can be asked replacing diameter by area or capacity.

Also in this paper we
considered analytic maps $f$ of the unit disk $\D$ into a region with 
bounded area, diameter or capacity, and established analogs of Schwarz's 
Lemma. What about analogs of Schwarz's Lemma for the `dual' situation 
of an analytic map $f:\Omega \to \D$, where $\Omega$ satisfies some 
geometric restriction?

\def\cprime{$'$}

\end{document}